\input amstex
\documentstyle{amsppt}
\input xypic

\def \ni{\noindent}

\def \ov{\overline}

\def \wt{\widetilde}
\def \ot{\otimes}

\def \sub{\subseteq}

\def \al{\alpha}

\def \de{\delta}
\def \De{\Delta}
\def \ep{\epsilon}

\def \si{\sigma}

\def \bx{\bold x}

\def \bN{\Bbb N}
\def \bZ{\Bbb Z}

\def \Alt{\operatorname{Alt}}
\def \Cyc{\operatorname{Cyc}}
\def \dim{\operatorname{dim}}
\def \Ext{\operatorname{Ext}}
\def \H{\operatorname{H}}
\def \HH{\operatorname{HH}}
\def \Hom{\operatorname{Hom}}

\def \ord{\operatorname{ord}}
\def \Tor{\operatorname{Tor}}

\def \ch{\operatorname{char}}

\topmatter
\title
Cohomology of trivial extensions of Frobenius algebras
\endtitle

\author
       Jorge A. Guccione and Juan J. Guccione
\endauthor

\address
     Jorge Alberto Guccione, Departamento de Matem\'atica, Facultad de
     Ciencias Exactas y Naturales, Pabell\'on 1 - Ciudad Universitaria,
     (1428) Buenos Aires, Argentina.
\endaddress

\email
     vander\@dm.uba.ar
\endemail

\address
     Juan Jos\'e Guccione, Departamento de Matem\'atica, Facultad de
     Ciencias Exactas y Naturales, Pabell\'on 1 - Ciudad Universitaria,
     (1428) Buenos Aires, Argentina.
\endaddress

\email
     jjgucci\@dm.uba.ar
\endemail

\abstract
     We obtain a decomposition for the Hochschild cochain complex of a
     split algebra and we study some properties of the cohomology of each
     term of this decomposition. Then, we consider the case of trivial
     extensions, specially of Frobenius algebras. In particular, we
     determine completely the cohomology of the trivial extension of a
     finite dimensional Hopf algebra. Finally, as an application, we obtain
     a result about the Hochschild cohomology of Frobenius algebras.

\endabstract



\thanks
Supported by UBACYT 01/TW79 and CONICET
\endthanks

\endtopmatter

\document

\head Introduction \endhead

Let $k$ be a field, $A$ a $k$-algebra and $M$ an $A$-bimodule. The
split algebra $E = A\ltimes M$, of $A$ with $M$, is the direct sum
$A\oplus M$ with the associative algebra structure given by
$$
(a+m)(a'+m') = aa' + am' + ma'.
$$
The (co)homology of split algebras, and of several particular types of
split algebras, such as the triangular matrix algebras and the trivial
extensions, has been considered in several papers. See for instance
\cite{C}, \cite{C-M-R-S}, \cite{G-G1}, \cite{G-M-S}, \cite{H1}, \cite{M-P}
and \cite{Mi-P}. This study is motivated in part by the relations between
the degree one Hochschild cohomology and the representation theory of a
finite dimensional algebra \cite{H2}, \cite{M-P}, \cite{S}, the relations
between the second and third Hochschild cohomology groups and the theory of
deformations of algebras \cite{G}, \cite{G-S}, and the following question
of Happel \cite{H1}: if an algebra has only a finite number of no nulls
Hochschild cohomology groups, is the algebra of finite homological
dimension? Moreover, the Hochschild cohomology groups are interesting
invariants of an algebra in its self.

In this work we continue the study of the Hochschild cohomology of split
algebras, computing the Hochschild cohomology of a trivial extension of a
Frobenius algebra $A$ in terms of the Hochschild cohomology of $A$.

Let $E=A\ltimes M$ a split algebra. This paper is organized as follows: in
Section~1 we show that the canonical Hochschild complex of $E$ with
coefficients in $E$ has a canonical decomposition as a direct sum of
subcomplexes $X^*_{(p)}$, such that $X^n_{(0)}$ is the Hochschild cochain
complex of $A$ with coefficients in $M$ and $X^n_{(p)}=0$ for all $n<p-1$.
Hence,
$$
\HH^n(E) = \H^n(A,M)\oplus \bigoplus_{p=1}^{n+1} H^n(X^*_{(p)}).
$$
Moreover, we prove that, under suitable hypothesis, for all $p\ge 1$ there
is a long exact sequence relating the cohomology of $X^*_{(p)}$, the $\Ext$
groups $\Ext_{A^e}^*(M^{\ot_A p-1},A)$ and the $\Ext$ groups $\Ext_{A^e}^*(
M^{\ot_A p},M)$.

In Section~2 we consider the cohomology of the trivial extension $TA$ of an
algebra $A$, which is the split algebra obtained taking $M=DA$, where $DA$
is the dual vector space of $A$, endowed with the usual $A$-bimodule
structure. For these algebras we compute $H^n(X^*_{(1)})$, for all $n\ge
0$, and $H^{p-1}(X^*_{(p)})$, for each $p\ge 2$.

The results of Sections~1 and 2 are close to the ones obtained in
\cite{C-M-R-S}.

In Section~3 we study trivial extensions of Frobenius algebras. The main
result is Theorem~3.10, where we compute the Hochschild cohomology of the
trivial extension of a finite order Frobenius $k$-algebra (see
Definition~3.7), when the characteristic of $k$ does not divide the order
of $A$. In particular, this result applies to the trivial extension of a
finite dimensional Hopf algebra. Finally, as an application, we obtain a
result about the Hochschild cohomology of a finite order Frobenius
$k$-algebra.

\head 1. A decomposition of the cohomology of an split algebra \endhead

Let $k$ be a field, $A$ a $k$-algebra and $M$ an $A$-bimodule. Let $E =
A\ltimes M$ the split algebra of $A$ with $M$ As it is well known, the
Hochschild cohomology $\HH^*(E) = \H(E,E)$ is the homology of the cochain
complex
$$
0 @>>> E @>b^1>> \Hom_k(E,E) @>b^2>> \Hom_k(E^{\ot 2},E) @>b^3>>
\Hom_k(E^{\ot 3},E) @>b^4>> \dots,
$$
where $E^{\ot n}$ is the $n$-fold tensor product of $E$, $b^1(x)(y) = yx
-xy$ and for $n>1$,
$$
\align
b^n(f)(x_1\ot\cdots\ot x_n) & = x_1f(x_2\ot\cdots\ot x_n)\\
& + \sum_{i=1}^{n-1} (-1)^i f(x_1\ot\cdots\ot x_ix_{i+1}\ot\cdots\ot x_n)\\
& + (-1)^n f(x_1\ot\cdots\ot x_{n-1})x_n.
\endalign
$$
For $0\le p\le n$, let $B^n_p\sub E^{\ot n}$ be the vector subspace spanned
by $n$-tensors $x_1\ot\cdots \ot x_n$ such that exactly $p$ of the $x_i$'s
belong to $M$, while the other $x_i$'s belong to $A$ (note that $B_0^0 =
k$). To unify expressions we make the convention that $B^n_p = 0$, for
$p<0$ or $n<p$. For each $p\ge 0$ we let $X^*_{(p)}$ denote the subcomplex
of $(\Hom_k(E^{\ot *},E),b^*)$ defined by:
$$
X^n_{(p)}:= \Hom_k(B^n_{p-1},A)\oplus \Hom_k(B^n_p,M).
$$
It is immediate that $(\Hom_k(E^{\ot *},E),b^*) = \bigoplus_{p=0}^{\infty}
X^*_{(p)}$. Hence, we have proved the following result:

\proclaim{Theorem 1.1} It is holds that
$$
\HH^n(E) = \bigoplus_{p=0}^{\infty} H^n(X^*_{(p)}) = \H^n(A,M)\oplus
\bigoplus_{p=1}^{n+1} H^n(X^*_{(p)}).
$$
\endproclaim

We also have the following result:

\proclaim{Theorem 1.2} If $\Tor_i^A(M,M^{\ot_A j}) = 0$ for $i>0$ and
$0<j<p$, then there is a long exact sequence
$$
\align
0 & \to  H^{p-1}(X^*_{(p)}) @>>> \Ext_{A^e}^0(M^{\ot_A p-1},A) @>>>
\Ext_{A^e} ^0(M^{\ot_A p},M) @>>> \\
& @>>> H^p(X^*_{(p)}) @>>> \Ext_{A^e}^1(M^{\ot_A p-1},A) @>>> \Ext_{A^e}^1(
M^{\ot_A p},M) @>>> \dots,
\endalign
$$
where, as usual, $M^{\ot_A 0} = A$.
\endproclaim

To prove this theorem, we need to study the cochain complexes $X^*_{(p)}$,
for $p\ge 1$. Let $\pi_A\:E\to A$ and $\pi_M\:E\to M$ be the maps
$\pi_A(a+m)= a$ and $\pi_M(a+m)=m$, respectively. Note that $X^*_{(p)}$ is
the total complex of the double complex
$$
X^{*,*}_{(p)}:=\quad
\diagram
&\vdots &\vdots \\
&\Hom_k(B^{p+1}_{p-1},A)\uto_{b^{0,p+2}_p}\rto^{\de^{1,p+1}_p} &
\Hom_k(B^{p+2}_p,M) \uto_{b^{1,p+2}_p}\\
&\Hom_k(B^p_{p-1},A)\uto_{b^{0,p+1}_p}\rto^{\de^{1,p}_p} &
\Hom_k(B^{p+1}_p,M) \uto_{b^{1,p+1}_p}\\
\text{row $p-1$}\save\go+<34pt,0pt>\Drop{\Text{$\to$}}\restore &
\Hom_k(B^{p-1}_{p-1},A) \uto_{b^{0,p}_p}\rto^{\de^{1,p-1}_p} &
\Hom_k(B^p_p,M) \uto_{b^{1,p}_p},\\
&\text{column $0$}\save \go+<0pt,18pt>\Drop{\Text{$\uparrow$}}\restore\\
\enddiagram
$$
where
$$
\align
\vspace{1.5\jot}
b_p^{0,n}(f)(x_1\ot\cdots\ot x_n) & = \pi_A(x_1)f(x_2\ot\cdots\ot x_n)\\
& + \sum_{i=1}^{n-1} (-1)^i f(x_1\ot\cdots\ot x_ix_{i+1}\ot\cdots\ot x_n)\\
& + (-1)^n f(x_1\ot\cdots\ot x_{n-1})\pi_A(x_n), \\
\vspace{1.5\jot}
b_p^{1,n}(g)(x_1\ot\cdots\ot x_{n+1}) & = x_1 g(x_2\ot\cdots\ot x_{n+1})\\
& + \sum_{i=1}^n(-1)^i g(x_1\ot\cdots\ot x_ix_{i+1}\ot\cdots\ot x_{n+1})\\
& + (-1)^{n+1} g(x_1\ot\cdots\ot x_n) x_{n+1},\\
\vspace{1.5\jot}
\de_p^{1,n-1}(f)(x_1\ot\cdots\ot x_n) & = \pi_M(x_1)f(x_2\ot\cdots\ot x_n)\\
& + (-1)^n f(x_1\ot\cdots\ot x_{n-1})\pi_M(x_n),
\endalign
$$
for $f\in \Hom_k(B^{n-1}_{p-1},A)$ and $g\in \Hom_k(B^n_p,M)$.

\smallskip

It is immediate that the cohomology $H^*(X^{0,*}_{(1)})$ of the $0$-column
of $X^{*,*}_{(1)}$ is the Hochschild cohomology $\HH^*(A)$. Moreover, we
have the following result:

\proclaim{Theorem 1.3} The following assertions hold

\smallskip

\item{1)} If $\Tor_i^A(M,M^{\ot_A j}) = 0$ for $i>0$ and $0<j<p-1$, then
$$
H^n(X^{0,*}_{(p)}) = \Ext_{A^e}^{n-p+1}(M^{\ot_A p-1},A)\quad\text{for all
$n\ge p-1$},
$$
\smallskip

\item{2)} If $\Tor_i^A(M,M^{\ot_A j}) = 0$ for $i>0$ and $0<j<p$, then
$$
H^n(X^{1,*}_{(p)}) = \Ext_{A^e}^{n-p+1}(M^{\ot_A p},M)\quad\text{for all
$n\ge p-1$}.
$$
\endproclaim

\demo{Proof} We prove the second assertion. The first one follows
similarly. In the proof of Theorem~2.5 of \cite{C-M-R-S} was showed that,
under our hypothesis, the complex
$$
M^{\ot_A p} @ <\mu^p<< A\ot B^p_p\ot A @<{b'_1}^p<< A\ot B^{p+1}_p \ot A
@<{b'_2}^p<<A\ot B^{p+2}_p \ot A @<{b'_3}^p<< \dots,
$$
where $\mu^p(a\ot x_1\ot\cdots\ot x_p\ot a') = ax_1\ot_A x_2\ot_A \cdots
\ot_A x_{p-1}\ot_A x_pa'$ and
$$
\align
{b'_n}^p(x_0\ot\cdots\ot x_{n+p+1}) & = \pi_A(x_0x_1)\ot x_2\ot\cdots\ot
x_{n+p+1}\\
& + \sum_{i=1}^{n+p-1} (-1)^i x_0\ot\cdots\ot x_ix_{i+1}\ot\cdots\ot
x_{n+p+1} \\
& + (-1)^{n+p} x_0\ot\cdots\ot x_{n+p-1}\ot \pi_A(x_{n+p}x_{n+p+1}),
\endalign
$$
is a projective resolution of $M^{\ot_A p}$ as an $A$-bimodule. The
assertion follows from this fact, using that $X^{1,*}_{(p)} \simeq
\Hom_{A^e}((A\ot B^{*+1}_p\ot A,{b'_{*-p+1}}^p),M)$.\qed

\enddemo

\subhead Proof of Theorem~1.2\endsubhead It follows from the long exact
sequence of homology of the short exact sequence
$$
0\to X^{1,*-1}_{(p)} @>>> X^*_{(p)} @>>> X^{0,*}_{(p)} \to 0,
$$
using Theorem~1.3.\qed

\remark{Remark 1.4} The vector spaces $B_p^n$ were considered in
\cite{C-M-R-S} in order to organize the canonical Hochschild cochain
complex of $E$ as a double complex. Using the decomposition $(\Hom_k(E^{\ot
*},E),b^*) = \bigoplus_{p=0}^{\infty} X^*_{(p)}$ obtained above
Theorem~1.1, can be easily shown that the spectral sequence introduced at
the beginning of Section~3 of \cite{C-M-R-S} satisfies $E^2 = E^3 = \cdots
= E^{\infty}$. However, the term $E^2$ is hard to compute in general.
\endremark

\head 2. Trivial extensions \endhead

Given an $A$-bimodule $M$, we let $DM$ denote $\Hom_k(M,k)$ endowed with
the usual $A$-bimodule structure.

\definition{Definition 2.1} Let $A$ be a $k$-algebra. The trivial extension
$TA$ of $A$ is the split algebra $A\ltimes DA$.
\enddefinition

\proclaim{Theorem 2.2} For each trivial extension $TA$, it is hold that
$$
H^0(X^*_{(1)}) = \HH^0(A) \quad\text{and}\quad H^n(X^*_{(1)}) = \HH^n(A)
\oplus \Ext_{A^e}^{n-1}(DA,DA)\,\,\forall n\ge 1.
$$
\endproclaim

\demo{Proof} Let $\de_1^{1,*}\:(\Hom_k(B^*_0,A),-b^{0,*+1}_1)\to
(\Hom_k(B^{*+1}_1,DA),b^{1,*+1}_1)$ be the map defined by
$$
\align
\de_1^{1,n-1}(f)(x_1\ot\cdots\ot x_n) & = \pi_{DA}(x_1) f(x_2\ot\cdots\ot
x_n)\\
& + (-1)^n f(x_1\ot\cdots\ot x_{n-1})\pi_{DA}(x_n).
\endalign
$$
By Theorem~1.3, $H^n(X^{0,*}_{(1)}) = \HH^n(A)$ and $H^{n-1}(X^{1,*}_{(1)})
= \Ext_{A^e}^{n-1}(DA,DA)$. Thus, we must prove that $H^n(X^*_{(1)}) =
H^n(X^{0,*}_{(1)})\oplus H^{n-1}(X^{1,*}_{(1)})$. Since $X^*_{(1)}$ is the
mapping cone of $\de_1^{1,*}$, to made out this task it suffices to check
that $\de_1^{1,*}$ is null homotopic. Let $\si_*\:\Hom_k (B^*_0,A)\to
\Hom_k(B^*_1,DA)$ be the family of maps defined by
$$
\si_n(f)(\bx_{1,n})(a) = (-1)^{jn+1} x_j(f(\bx_{j+1,n}\ot a\ot
\bx_{1,j-1})) \quad\text{if $x_j\in DA$,}
$$
where, to abbreviate, we write $\bx_{h,l} = x_h\ot \cdots \ot x_l$, for
$h<l$. We assert that $\si_*$ is an homotopy from $\de_1^{1,*}$ to $0$. We
have,
$$
\align
b^{1,n}_1(\si_n(f))(\bx_{1,n+1}) & = \pi_A(x_1)(\si_n(f)(\bx_{2,n+1}))
+ (-1)^{n+1} (\si_n(f)(\bx_{1,n}))\pi_A(x_{n+1})\\
& + \sum_{i=1}^n (-1)^i \si_n(f)(\bx_{1,i-1}\ot x_ix_{i+1}\ot \bx_{i+2,n+1}).
\endalign
$$
Hence, if $x_1\in DA$, then
$$
\align
b^{1,n}_1(\si_n(f))(\bx_{1,n+1})(x_{n+2}) & = (-1)^{n+2}x_1(x_2f(
\bx_{3,n+2})) \\
& + \sum_{i=2}^{n+1} (-1)^{n+i+1} x_1(f(\bx_{2,i-1}\ot x_ix_{i+1}\ot
\bx_{i+2,n+2})),
\endalign
$$
if $x_j\in DA$ for $1<j\le n$, then
$$
\align
b^{1,n}_1(\si_n(f))(&\bx_{1,n+1})(x_0) = (-1)^{(j-1)n+j} x_j(f(\bx_{j+1,
n+1} \ot \bx_{0,j-2})x_{j-1})\\
& + \sum_{i=0}^{j-2} (-1)^{(j-1)n+i +1} x_j(f(\bx_{j+1,n+1}\ot
\bx_{0,i-1}\ot x_ix_{i+1}\ot \bx_{i+2,j-1})) \\
& + (-1)^{jn+j+1}x_j(x_{j+1}f(\bx_{j+2,n+1}\ot\bx_{0,j-1}))\\
& + \sum_{i=j+1}^n (-1)^{jn+i+1} x_j(f(\bx_{j+1,i-1}\ot x_ix_{i+1}\ot
\bx_{i+2,n+1}\ot \bx_{0,j-1})) \\
& + (-1)^{jn+n} x_j(f(\bx_{j+1,n}\ot x_{n+1}x_0\ot \bx_{1,j-1})),
\endalign
$$
and if $x_{n+1}\in DA$, then
$$
\align
b^{1,n}_1(\si_n(f))(\bx_{1,n+1})(x_0) & = \sum_{i=0}^{n-1} (-1)^{n+i+1}
x_{n+1}(f(\bx_{0,i-1}\ot x_ix_{i+1}\ot \bx_{i+2,n})) \\
& - x_{n+1}(f(\bx_{0,n-1})x_n).
\endalign
$$
On the other hand, if $x_1\in DA$, then
$$
\align
\si_{n+1}(-b^{0,n}_1(f))(\bx_{1,n+1})(x_{n+2}) & = (-1)^{n+1}
x_1(b^{0,n}_1(f)(\bx_{2,n+2}))\\
& = (-1)^{n+1} x_1(x_2f(\bx_{3,n+2})) + x_1(f(\bx_{2,n+1})x_{n+2}) \\
& + \sum_{i=2}^{n+1} (-1)^{n+i} x_1(f(\bx_{2,i-1}\ot x_ix_{i+1}\ot
\bx_{i+2,n+2})),
\endalign
$$
if $x_j\in DA$ for $1<j\le n$, then
$$
\align
\si_{n+1}(-b^{0,n}_1&(f))(\bx_{1,n+1})(x_0) = (-1)^{j(n+1)+1}
x_1(b^{0,n}_1(f)(\bx_{j+1,n+1}\ot \bx_{0,j-1}\ot ))\\
& = (-1)^{j(n+1)}x_j(x_{j+1}f(\bx_{j+2,n+1}\ot\bx_{0,j-1}))\\
& + \sum_{i=j+1}^n (-1)^{j(n+1)+i-j} x_j(f(\bx_{j+1,i-1}\ot x_ix_{i+1}\ot
\bx_{i+2,n+1}\ot \bx_{0,j-1})) \\
& + (-1)^{j(n+1)+n-j+1} x_j(f(\bx_{j+1,n}\ot x_{n+1}x_0\ot \bx_{1,j-1}))\\
& + \sum_{i=0}^{j-2} (-1)^{j(n+1)+i+n-j} x_j(f(\bx_{j+1,n+1}\ot
\bx_{0,i-1}\ot x_ix_{i+1}\ot \bx_{i+2,j-1})) \\
& + (-1)^{j(n+1)+n+1} x_j(f(\bx_{j+1,n+1}\ot \bx_{0,j-2})x_{j-1}),
\endalign
$$
and if $x_{n+1}\in DA$, then
$$
\align
\si_{n+1}(-b^{0,n}_1(f))(\bx_{1,n+1})(x_0) & = (-1)^{n+1}
x_{n+1}(b^{0,n}_1(f)(\bx_{0,n}))\\
& = (-1)^{n+1} x_{n+1}(x_0f(\bx_{1,n})) + x_{n+1}(f(\bx_{0,n-1})x_n) \\
& + \sum_{i=0}^{n-1} (-1)^{n+i} x_{n+1}(f(\bx_{0,i-1}\ot x_ix_{i+1}\ot
\bx_{i+2,n})).
\endalign
$$
The assertion follows immediately from these equalities.\qed
\enddemo

As usual, for each $A$-bimodule $M$ we write $M\ot_{A^e} = M/[A,M]$, where
$[A,M]$ is the vector subspace of $M$ generated by $\{am - ma: a\in A,m\in
M\}$. The map $\theta\:\Hom_{A^e}((DA)^{\ot_A p-1},A) \to \Hom_k((DA)^{
\ot_A p}\ot_{A^e},k)$ defined by
$$
\theta(f)(\psi_1\ot\cdots\ot \psi_p) = \psi_p(f(\psi_1\ot\cdots\ot
\psi_{p-1}))
$$
is injective. For each $p\ge 2$, we let $\Cyc^p_A(DA)$ denote the set of
$k$-linear maps $g\:(DA)^{\ot_A p}\to k$ verifying $g(\psi_1\ot\cdots
\ot\psi_p) = (-1)^{p-1} g(\psi_2\ot\cdots \ot\psi_p\ot\psi_1)$. Note that
$\Cyc^p_A(DA) \sub \Hom_k((DA)^{\ot_A}\ot_{A^e},k)$.

\smallskip

The following lemma and its proof is inspired in the proof of Theorem~5.5
of \cite{C-M-R-S}.

\proclaim{Lemma 2.3} For each $p\ge 2$, it holds that
$$
H^{p-1}(X^*_{(p)})\simeq \theta(\Hom_{A^e}((DA)^{\ot_A p-1},A))\cap
\Cyc^p_A(DA).
$$
Moreover, if $A$ is a finite dimensional $k$-algebra, then $\theta(
\Hom_{A^e}((DA)^{\ot_A p-1},A))\cap \Cyc^p_A(DA) = \Cyc^p_A(DA)$.

\endproclaim

\demo{Proof} Let us compute $H^{p-1}(X^*_{(p)})$. It is easy to see that
$H^{p-1}(X^*_{(p)})$ is the kernel of
$$
\wt{\de}_p^{1,p-1}\: \Hom_{A^e}((DA)^{\ot_A p-1},A)\to \Hom_{A^e}((DA)^{
\ot_A p},DA),
$$
where $\wt{\de}_p^{1,p-1}$ is the map defined by $\wt{\de}_p^{1,p-1}(f)
(\psi_1\ot\cdots \ot\psi_p) = \psi_1 f(\psi_2\ot\cdots \ot\psi_p) +
(-1)^p f(\psi_1\ot\cdots \ot\psi_{p-1})\psi_p$. Consider the isomorphism
$$
\vartheta\:\Hom_{A^e}((DA)^{\ot_A p},DA) \to \Hom_k((DA)^{\ot_A
p}\ot_{A^e},k),
$$
given by $\vartheta(f)(\psi_1\ot\cdots \ot\psi_p) = f(\psi_1\ot\cdots
\ot\psi_p)(1)$. Let $g$ be in the the image of $\theta$. We have
$$
\align
&(\vartheta\circ \wt{\de}_p^{1,p-1}\circ \theta^{-1}(g))(\psi_1\ot\cdots\ot
\psi_p) \\
& = \psi_1 (\theta^{-1}(g)(\psi_2\ot\cdots \ot\psi_p)) + (-1)^p
\psi_p(\theta^{-1}(g)(\psi_1\ot\cdots \ot\psi_{p-1}))\\
& = g(\psi_2\ot\cdots \ot\psi_p\ot\psi_1) + (-1)^p g(\psi_1\ot\cdots
\ot\psi_p).
\endalign
$$
Hence $H^{p-1}(X^*_{(p)}) \simeq \theta(\Hom_{A^e}((DA)^{\ot_A p-1},A))\cap
\Cyc^p_A(DA)$, as desired. To finish the proof it is sufficient to note
that if $\dim A <\infty$, then  $\Hom_k((DA)^{\ot_A p}\ot_{A^e},k)\simeq
\Hom_{A^e}((DA)^{\ot_A p-1},DDA) \simeq \Hom_{A^e}((DA)^{\ot_A p-1},A)$,
which implies that $\theta$ is an isomorphism.\qed

\enddemo

The following result improves Theorem~5.7 of \cite{C-M-R-S}.

\proclaim{Corollary 2.4} For all algebra $A$ and each $n\ge 1$, we have
$$
\align
\HH^n(TA) & = \HH^n(A)\oplus \HH_n(A)^*\oplus \Ext_{A^e}^{n-1}(DA,DA)\\
& \oplus \theta(\Hom_{A^e}((DA)^{\ot_A n},A))\cap \Cyc^{n+1}_A(DA) \oplus
\bigoplus_{p=2}^n H^n(X^*_{(p)}).
\endalign
$$
\endproclaim

\demo{Proof} It follows from Theorem~1.1, Theorem~2.2, Lemma~2.3 and the
fact that $\H^n(A,DA) = \HH_n(A)^*$.\qed
\enddemo

\proclaim{Corollary 2.5} Let $A$ be a finite dimensional $k$-algebra. For
each $n\ge 1$, we have
$$
\HH^n(TA) = \HH^n(A)\oplus \HH_n(A)^*\oplus \Ext_{A^e}^{n-1}(DA,DA)
\oplus \Cyc^{n+1}_A(DA) \oplus \bigoplus_{p=2}^n H^n(X^*_{(p)}).
$$
\endproclaim

\demo{Proof} It follows from Corollary~2.4 and Lemma~2.3.\qed
\enddemo

\proclaim{Lemma 2.6} Let $(DDA)^A = \{\varphi\in DDA : \text{ $a\varphi =
\varphi a$ for all $a\in A$}\}$. It is hold that $\Ext^0_{A^e}(DA,DA) =
(DDA)^A$.
\endproclaim

\demo{Proof} We have
$$
\align
\Ext^0_{A^e}(DA,DA) & = \Hom_{A^e}(DA,DA) \simeq \Hom_k(DA\ot_{A^e} A,k)\\
& \simeq \Hom_k(A\ot_{A^e} DA,k)\simeq \Hom_{A^e}(A,DDA) = (DDA)^A.\qed
\endalign
$$
\enddemo

The following two results were obtained in \cite{C-M-R-S}. They use the
notation $\Alt_A(DA)$ instead of $\Cyc^2_A(DA)$.

\proclaim{Theorem 2.7} For each trivial extension $TA$, it is hold that:

\smallskip

\item{1)} $\HH^0(TA) = \HH^0(A) \oplus \HH_0(A)^*$,

\smallskip

\item{2)} $\HH^1(TA) = \HH^1(A) \oplus \HH_1(A)^*\oplus (DDA)^A\oplus
\theta(\Hom_{A^e}(DA,A))\cap \Cyc^2_A(DA)$.

\endproclaim

\demo{Proof} 1) By Theorems~1.1 and 2.2, and the fact that $\H^n(A,DA) =
\HH_n(A)^*$, we have
$$
\HH^0(TA) = H^0(A,DA) \oplus H^0(X^*_{(1)}) = \HH_0(A)^* \oplus \HH^0(A)
$$
\smallskip

\ni 2) It follows immediately from Corollary~2.4 and Lemma~2.6.\qed

\enddemo

\proclaim{Corollary 2.8} Let $TA$ be a trivial extension of a finite
dimensional $k$-algebra. Then,

\smallskip

\item{1)} $\HH^0(TA) = \HH^0(A) \oplus \HH_0(A)^*$,

\smallskip

\item{2)} $\HH^1(TA) = \HH^1(A) \oplus \HH_1(A)^*\oplus A^A\oplus
\Cyc^2_A(DA)$.
\endproclaim

\demo{Proof} It follows immediately from Theorem~2.7, Lemma~2.3 and the
fact that $DAA \simeq A$.\qed
\enddemo

\head 3. Trivial extensions of Frobenius algebras \endhead

We recall that a finite dimensional $k$-algebra is Frobenius if there
exists a linear form $\varphi\: A\to k$ such that the map $A\to DA$,
defined by $x\mapsto x\varphi$ is a left $A$-module isomorphism. This
linear form $\varphi\: A\to k$ is called a Frobenius homomorphism. It is
well known that this is equivalent to say that the map $x\mapsto \varphi
x$, from $A$ to $DA$, is an isomorphism of right $A$-modules. From this
follows easily that there exists an automorphism $\rho$ of $A$, called the
Nakayama automorphism of $A$ with respect to $\varphi$, such that $x\varphi
= \varphi \rho(x)$, for all $x\in A$. It is easy to check that a linear
form $\wt{\varphi}\: A\to k$ is another Frobenius homomorphism if and only
if there exists $x\in A$ invertible, such that $\wt{\varphi} = x\varphi$.
Also it is easy to check that the Nakayama automorphism of $A$ with respect
to $\wt{\varphi}$ is the map given by $a\mapsto \rho(x)^{-1}\rho(a)\rho(x)$.

Given algebra maps $f$ and $g$, we let $A_f^g$ denote $A$ endowed with
the $A$-bimodule structure given by $a\cdot x\cdot b : = f(a)xg(b)$. To
simplify notations we write $A_f$ instead of $A_f^{id}$ and $A^g$ instead
of $A_{id}^g$. We have the $A$-bimodule isomorphism $\Theta\:(DA)^{\ot_A p}
\to A_{\rho^p}$, given by $\Theta(\varphi x_1\ot_A\cdots \ot_A \varphi x_p)
= \rho^{p-1} (x_1) \rho^{p-2}(x_2)\cdots \rho(x_{p-1})x_p$. Let
$$
A_{\rho^p} @ <\mu<< A\ot A_{\rho^p}  @<b'_1<< A^{\ot 2}\ot A_{\rho^p}
@<b'_2<< A^{\ot 3}\ot A_{\rho^p} @<b'_3<< A^{\ot 4}\ot A_{\rho^p} @<b'_4<<
\dots,
$$
be the bar resolution of $A_{\rho^p}$.

\proclaim{Proposition 3.1} Let $(A\ot B^{*+p}_p\ot A,{b'_*}^p)$ be as in
the proof of Theorem~1.3. The following facts hold:

\smallskip

\item{1)} There is a chain map $\Theta^p_{*+p}\: (A\ot B^{*+p}_p\ot
A,{b'_*}^p)\to (A^{\ot *+1}\ot A_{\rho^p},b'_*)$, given by
$$
\Theta^p_{n+p}(\bx_{0,n+p+1}) = \cases\!\bx_{0n} \ot \Theta(x_{n+1}\ot_A
\cdots \ot_A x_{n+p})x_{n+p+1} & \text{if $x_1,\dots,x_n \in A$,} \\ \! 0
&\text{in other case,} \endcases
$$
where $\bx_{0,n+p+1} = x_0\ot \cdots \ot x_{n+p+1}$ and $\bx_{0n} =
x_0\ot\cdots\ot x_n$,

\smallskip

\item{2)} There is a chain map $\Psi^p_{*+p}\: (A^{\ot *+1}\ot
A_{\rho^p},b'_*)\to (A\ot B^{*+p}_p\ot A,{b'_*}^p)$, given by
$$
\align
\Psi^p_{n+p}&(x_0\ot \cdots \ot x_{n+1}) = \sum_{0\le i_1\le \cdots \le
i_p\le n} (-1)^{i_1+\cdots+i_p+pn} x_0\ot\cdots \ot x_{i_1}\ot\varphi\\
&\ot\rho(x_{i_1+1})\ot\cdots\ot\rho(x_{i_2})\ot\varphi\ot \rho^2(x_{i_2+1})
\ot\cdots \ot \rho^2(x_{i_3}) \ot\varphi\ot\cdots\ot\varphi\\
& \ot\rho^{p-1}(x_{i_{p-1}+1})\ot\cdots \ot \rho^{p-1}(x_{i_p})\ot\varphi
\ot \rho^p(x_{i_p+1})\ot\cdots\ot\rho^p(x_n)\ot x_{n+1}.
\endalign
$$

\smallskip

\item{3)} $\mu\circ \Theta^p_p = \Theta\circ \mu^p$ and $\Theta\circ \mu^p
\circ \Psi^p_p = \mu$, where $\mu^p\:A\ot B^p_p\ot A\to (DA)^{\ot_A p}$ is
the map introduce in the proof of Theorem~1.3.

\endproclaim

\demo{Proof} We left items 1) and 3) to the reader. Let us see 2). For $0\le
i_1< \cdots < i_p \le n$ we write
$$
\align
T&_{i_1,\dots,i_p} = (-1)^{i_1+\cdots+i_p+pn} x_0\ot\cdots \ot x_{i_1}\ot
\varphi \ot\rho(x_{i_1+1})\ot\cdots\ot\rho(x_{i_2})\ot \varphi\ot\cdots\\
& \ot\varphi\ot\rho^{p-1}(x_{i_{p-1}+1})\ot\cdots\ot \rho^{p-1}(x_{i_p})\ot
\varphi \ot \rho^p(x_{i_p+1})\ot\cdots\ot\rho^p(x_n)\ot x_{n+1}.
\endalign
$$
The term of $b'_n\circ \Psi^p_{n+p}(x_0\ot \cdots \ot x_{n+1})$ obtained
multiplying $\rho^{j-1}(x_{i_j})$ by $\varphi$ in $T_{i_1,\dots,i_p}$ (with
$i_j>0$) cancels with the term of $b'_n\circ \Psi^p_{*+p}(x_0\ot \cdots \ot
x_{n+1})$ obtained multiplying $\varphi$ by $\rho^j(x_{i_j})$ in
$T_{i_1,\dots,i_{j-1},i_j-1,i_{j+1},\dots,i_p}$. Using this fact it is easy
to see that ${b'_n}^p\circ \Psi^p_{n+p}(x_0\ot \cdots \ot x_{n+1}) =
\Psi^p_{n+p-1}\circ b'_n(x_0\ot \cdots \ot x_{n+1})$.\qed
\enddemo

\proclaim{Theorem 3.2} Let $p\ge 1$. Assume that $A$ is a Frobenius
algebra. Then, for the trivial extension $TA$, the double complex
$X^{*,*}_{(p)}$, introduced below Theorem~1.2, have the same homology
groups that the complex
$$
Y^{*,*}_{(p)}\:\quad
\diagram
&\vdots &\vdots \\
&\Hom_k(A^{\ot 2},A)\uto_{\wt{b}^{0,p+2}_p} \rto^{\wt{\de}^{1,p+1}_p} &
\Hom_k(A^{\ot 2},A) \uto_{\wt{b}^{1,p+2}_p}\\
&\Hom_k(A,A)\uto_{\wt{b}^{0,p+1}_p} \rto^{\wt{\de}^{1,p}_p} & \Hom_k(A,A)
\uto_{\wt{b}^{1,p+1}_p}\\
\text{row $p-1$}\save\go+<34pt,0pt>\Drop{\Text{$\to$}}\restore &
\Hom_k(k,A) \uto_{\wt{b}^{0,p}_p}\rto^{\wt{\de}^{1,p-1}_p} & \Hom_k(k,A)
\uto_{\wt{b}^{1,p}_p},\\
&\text{column $0$}\save \go+<0pt,18pt>\Drop{\Text{$\uparrow$}}\restore\\
\enddiagram
$$
where
$$
\allowdisplaybreaks
\align
\vspace{1.5\jot}
\wt{b}_p^{0,n}(f)(x_1\ot\cdots\ot x_{n-p+1}) & = x_1f(x_2\ot\cdots \ot
x_{n-p+2})\\
& + \sum_{i=1}^{n-p} (-1)^i f(x_1\ot\cdots\ot x_ix_{i+1}\ot\cdots\ot
x_{n-p+1})\\
& + (-1)^{n-p+1} f(x_1\ot\cdots\ot x_{n-p}) \rho^{p-1}(x_{n-p+1}),\\
\vspace{1.5\jot}
\wt{b}_p^{1,n}(g)(x_1\ot\cdots\ot x_{n-p+1}) & = x_1 g(x_2\ot\cdots\ot
x_{n-p+1})\\
& + \sum_{i=1}^{n-p}(-1)^i g(x_1\ot\cdots\ot x_ix_{i+1}\ot\cdots\ot
x_{n-p+1})\\
& + (-1)^{n-p+1} g(x_1\ot\cdots\ot x_{n-p}) \rho^{p-1}(x_{n-p+1}),\\
\vspace{1.5\jot}
\wt{\de}_p^{1,n-1}(f)(x_1\ot\cdots\ot x_{n-p}) & = (-1)^n
f(x_1\ot\cdots\ot x_{n-p})\\
& + (-1)^{n-p} \rho^{-1}(f(\rho(x_1)\ot\cdots\ot \rho(x_{n-p}))),
\endalign
$$
for $f,g\in \Hom_k(A^{\ot n-p},A)$.

\endproclaim

\demo{Proof} Let $\Theta^{p-1}_*$ and $\Psi^p_*$ be as in Proposition~3.1.
From Theorem~1.3 and the fact that the bar resolutions of $A_{\rho^{p-1}}$
and $A_{\rho^p}$ are $A^e$-projective resolutions, it follows that the maps
$$
\Theta_{(p)}^{0,*} = \Hom_{A^e}(\Theta^{p-1}_*,A)\quad\text{and}\quad
\Psi_{(p)}^{1,*} = \Hom_{A^e}(\Psi^p_*,DA)
$$
are quasiisomorphism from $\Hom_{A^e}\bigl((A^{\ot *-p+2}\ot
A_{\rho^{p-1}},b'_*),A\bigr)$ to $X_{(p)}^{0,*}$ and from $X_{(p)}^{1,*}$
to $\Hom_{A^e}\bigl((A^{\ot *-p+2}\ot A_{\rho^p},b'_*),DA\bigr)$
respectively. Let
$$
\align
& \Xi_{(p)}^{0,*} \: Y^{0,*}_{(p)} @>>> \Hom_{A^e}\bigl((A^{\ot *-p+2}\ot
A_{\rho^{p-1}},b'_*),A\bigr)\\
\intertext{and}
& \Upsilon_{(p)}^{1,*}\: Y^{1,*}_{(p)} @>>> \Hom_{A^e}\bigl((A^{\ot *-p+2}
\ot A_{\rho^p},b'_*),DA\bigr),
\endalign
$$
be the chain complex isomorphisms defined by
$$
\align
& \Xi_{(p)}^{0,n}(f)(a_0\ot\cdots\ot a_{n-p+2}) = a_0f(a_1\ot\cdots\ot
a_{n-p+1})a_{n-p+2}
\intertext{and}
& \Upsilon_{(p)}^{1,n}(f)(a_0\ot\cdots\ot a_{n-p+2}) = a_0\varphi
\rho(f(a_1\ot\cdots\ot a_{n-p+1}))a_{n-p+2}.
\endalign
$$
To finish the proof it suffices to check that $\Upsilon_{(p)}^{1,*}\circ
\wt{\de}_p^{1,*} = \Psi_{(p)}^{1,*}\circ\de_p^{1,*}\circ \Theta_{(p)}^{0,*}
\circ \Xi_{(p)}^{0,*}$.\qed
\enddemo

\remark{Remark 3.3} Recall that $A^{\rho^{p-1}}$ is $A$ endowed with the
$A$-bimodule structure given by $a\cdot x\cdot b = ax\rho^{p-1}(b)$. The
columns of $Y^{*,*}_{(p)}$ are the canonical Hochschild cochain complex of
$A$ with coefficients in $A^{\rho^{p-1}}$. Hence, $H^n(Y^{0,*}_{(p)}) =
H^n(Y^{1,*}_{(p)}) = \H^n(A,A^{\rho^{p-1}})$.
\endremark

\medskip

For each $p\ge 1$ we let $Y^*_{(p)}$ denote the total complex of
$Y^{*,*}_{(p)}$.

\proclaim{Proposition 3.4} Assume that $A$ is a Frobenius algebra. Then,
for the trivial extension $TA$, it is hold that
$$
H^n(Y^*_{(1)}) = H^n(X^*_{(1)}) = \cases \HH^n(A) & \text{if $n=0$,}\\
\HH^n(A) \oplus \HH^{n-1}(A) & \text{if $n>0$.}\endcases
$$
\endproclaim

\demo{Proof} By Theorems~3.2, 1.3 and the proof of Theorem~2.2, we get
$$
H^n(Y^*_{(1)}) = H^n(X^*_{(1)}) = H^n(X^{0,*}_{(1)})\oplus
H^{n-1}(X^{1,*}_{(1)}) = \HH^n(A) \oplus H^{n-1}(X^{1,*}_{(1)}).
$$
Moreover, by the proof of Theorem~3.2, we know that $Y^{1,*}_{(1)} \simeq
X^{1,*}_{(1)}$. Now, the result follows from Remark~3.3.\qed
\enddemo

\proclaim{Corollary 3.5} Assume that $A$ is a Frobenius algebra. Then, for
all $n\ge 1$, we have
$$
\HH^n(TA) = \HH_n(A)^* \oplus\HH^n(A) \oplus \HH^{n-1}(A)\oplus
\Cyc^{n+1}_A(DA) \oplus \bigoplus_{p=2}^n H^n(Y^*_{(p)}).
$$
\endproclaim

\demo{Proof} It follows immediately from Theorem~1.1, Lemma~2.3,
Proposition~3.4, Theorem~3.2 and the fact that $\H^n(A,DA) =
\HH_n(A)^*$.\qed
\enddemo

\remark{Remark 3.6} Let $A$ be a Frobenius algebra. By Lemma~2.3 and
Theorem~3.2,
$$
\align
\Cyc^n_A(DA) & = H^n(X^*_{(n+1)}) = H^n(Y^*_{(n+1)})\\
& = \{x\in A : \rho(x) = (-1)^n x \text{ and $ax = x\rho^n(a)$ for all
$a\in A$}\}.
\endalign
$$
\endremark

\definition{Definition 3.7} Let $A$ be a Frobenius $k$-algebra, $\varphi\:
A\to k$ a Frobenius homomorphism and $\rho\:A\to A$ the Nakayama
automorphism with respect to $\varphi$. We say that $A$ has order $m\in\bN$
and we write $\ord_A=m$ if $\rho^m= id_A$ and $\rho^r\ne id_A$ for all
$r<m$. In this case we also write
$$
e_A=\cases \ord_A & \text{if $\ord_A$ is even,}\\
2\ord_A & \text{if $\ord_A$ is odd.}\endcases
$$
If $\rho^m\ne id_A$ for all $m\in \bN$, then we say that $A$ has infinite
order. By the discussion at the beginning of this section these definitions
are independent of $\varphi$.
\enddefinition

From now on $A$ is a Frobenius algebra of finite order.

\proclaim{Theorem 3.8} For each $n\ge 1$ we write $n = qe_A+s$, with $0\le
s< e_A$. Then,
$$
\align
\HH^n(TA) & = \HH_n(A)^* \oplus\bigoplus_{i=0}^q
\left(\HH^{n-ie_A}(A)\oplus \HH^{n-ie_A-1}(A)\right) \\
& \oplus\bigoplus_{j=2}^{e_A}\bigoplus_{i=0}^{q-1} H^{n-ie_A} (Y^*_{(j)})
\oplus\bigoplus_{j=2}^{s+1} H^{n-qe_A}(Y^*_{(j)}),
\endalign
$$
where $\HH^{-1}(A) = 0$.

\endproclaim

\demo{Proof} Since $\rho^{e_A} = id$ and $e_A$ is even, we have
$Y^{*,*}_{(p)}=Y^{*,*-e_A}_{(p-e_A)}$ for all $p\ge e_A+1$. Hence, if $1\le
r\le e_A$ is such that $p \equiv r \pmod{e_A}$, then $H^n(Y^*_{(p)}) =
H^{n-p+r}(Y^*_{(r)})$ for all $p\ge 1$. By Theorem~1.1, the fact that
$\H^n(A,DA) = \HH_n(A)^*$, Theorem~3.2 and the above equality, we have
$$
\align
& \HH^n(TA) = \HH_n(A)^* \oplus \bigoplus_{p=1}^{n+1} H^n(Y^*_{(p)}) \\
& = \HH_n(A)^* \oplus \bigoplus_{j=1}^{s+1} \bigoplus_{i=0}^q
H^{n-ie_A}(Y^*_{(j)}) \oplus\bigoplus_{j=s+2}^{e_A} \bigoplus_{i=0}^{q-1}
H^{n-ie_A}(Y^*_{(j)})\\
& = \HH_n(A)^*\oplus H^{n-qe_A}(Y^*_{(s+1)}) \oplus \bigoplus_{j=1}^{e_A}
\bigoplus_{i=0}^{q-1} H^{n-ie_A}(Y^*_{(j)})\oplus\bigoplus_{j=1}^s
H^{n-qe_A}(Y^*_{(j)}).
\endalign
$$
The assertion follows immediately from Proposition~3.4.\qed

\enddemo

Given  an $A^e$-projective resolution $(X_*,\partial_*)$ of $A$, let $X'_n
= X_n$ endowed with the $A$-bimodule action $a\cdot x \cdot b :=
\rho(a)x\rho(b)$. Then $(X'_*,\partial_*)$ is an $A^e$-resolution of
$A^{\rho}_{\rho}$. Since $\rho\:A\to A^{\rho}_{\rho}$ is an $A$-bimodule
map, there is a chain map $\wt{\rho}_* \: (X_*,\partial_*) \to
(X'_*,\partial_*)$, unique up to chain homotopy equivalence, lifting
$\rho$.  Let $p\ge 1$. It is easy to check that the map
$$
f\mapsto \rho^{-1}\circ f\circ \wt{\rho}_n \qquad \text{$n\ge 0$ and $f\in
\Hom_{A^e}(X_n,A^{\rho^{p-1}})$,}
$$
is an automorphism of $\Hom_{A^e}\bigl((X_*,\partial_*), A^{\rho^{p-1}}
\bigr)$. Taking cohomology we obtain an automorphism $\wt{\rho}^*$ of
$\H^*(A,A^{\rho^{p-1}})$, which do not depend of the choose resolution.
Clearly $(\wt{\rho}^*)^{\ord_A}$ is the identity map. Assume that the
characteristic of $k$ does not divide $\ord_A$ and that $k$ has a
primitive $\ord_A$-th root of unity $w$. Since $X^{\ord_A} - 1$ has
distinct roots $w^i$ ($0\le i<\ord_A$), the cohomology
$\H^n(A,A^{\rho^{p-1}})$ decomposes as the direct sum
$$
\H^n(A,A^{\rho^{p-1}}) = \bigoplus_{l=0}^{\ord_A-1}\H_{(p)}^{n,l}(A),
$$
where  $\H_{(p)}^{n,l}(A) = \{x\in \H^n(A,A^{\rho^{p-1}}):
\wt{\rho}^n(x) = w^{-l}x\}$.

\proclaim{Proposition 3.9} Under the above hypothesis, we have that for
each $n\ge p-1$,
$$
H^n(Y_{(p)}^{*,*}) = \! \cases \!\H_{(p)}^{n-p+1,0}(A) \oplus
\H_{(p)}^{n-p,0}(A) &\!\text{if $p$ is odd,}  \\
\!\H_{(p)}^{n-p+1,0}(A)\oplus \H_{(p)}^{n-p,0}(A) &\!\text{if $\ch(k)=2$
and $p$ is even,} \\
\!\H_{(p)}^{n-p+1,\frac{e_A}{2}}(A)\oplus \H_{(p)}^{n-p,\frac{e_A}{2}}(A)
&\!\text{if $\ch(k)\ne 2$ and $p,\ord_A$ are even,}\\
\! 0 &\!\text{in other case,} \endcases
$$
where $\HH_{(p)}^{-1,0}(A) = \HH_{(p)}^{-1,\frac{e_A}{2}}(A) = 0$.

\endproclaim

\demo{Proof} Since the minimal polynomial $X^{\ord_A} -1$ of $\rho\:A\to
A$, has distinct roots $w^i$ ($0\le i<\ord_A$), the algebra $A$ becomes a
$\frac{\bZ}{\ord_A\bZ}$-graded algebra
$$
A = A_0\oplus \cdots\oplus A_{\ord_A-1}, \quad\text{where $A_u = \{a\in
A:\rho(a) = w^u a\}$}.
$$
For each $0\le l < \ord_A$, let $Y_{(p),l}^{*,*}$ be the subcomplex of
$Y^{*,*}_{(p)}$ defined by
$$
Y_{(p),l}^{i,n} = \bigoplus_{B_{l,n}} \Hom(A_{u_1}\ot\cdots\ot
A_{u_{n-p+1}},A_v),
$$
where $B_{l,n} = \{(u_1,\dots,u_{n-p+1},v)$ such that $v - u_1 -\cdots
-u_{n-p+1} \equiv l \pmod{\ord_A}\}$. It is clear that $Y^{*,*}_{(p)} =
\bigoplus_{l=0}^{\ord_A} Y_{(p),l}^{*,*}$. Let $f\in Y_{(p),l}^{0,n}$. A
direct computation shows that
$$
\wt{\de}_p^{1,n}(f)(x_1\ot\cdots\ot x_{n-p+1}) = (-1)^{n+1}(1 + (-1)^p
w^{-l}) f(x_1\ot\cdots\ot x_{n-p+1}).
$$
Hence the horizontal boundary maps of $Y_{(p),l}^{*,*}$ are isomorphisms if
$w^l \ne (-1)^{p-1}$, and they are zero maps if $w^l = (-1)^{p-1}$. So,
$$
H^n(Y_{(p),l}^{*,*}) = \cases 0 &\text{if $w^l \ne (-1)^{p-1}$,}\\
H^n(Y_{(p),l}^{0,*}) \oplus H^{n-1}(Y_{(p),l}^{0,*}) &\text{if $w^l =
(-1)^{p-1}$.} \endcases
$$
Then,
$$
H^n(Y_{(p)}^{*,*}) = \! \cases \!H^n(Y_{(p),0}^{0,*}) \oplus
H^{n-1}(Y_{(p),0}^{0,*}) &\text{if $p$ is odd,}  \\
\!H^n(Y_{(p),0}^{0,*}) \oplus H^{n-1}(Y_{(p),0}^{0,*}) &\text{if
$\ch(k)=2$ and $p$ is even,} \\
\!H^n(Y_{(p),\frac{e_A}{2}}^{0,*})\oplus H^{n-1}(Y_{(p),\frac{e_A}{2}}^{0,*})
&\text{if $\ch(k)\ne 2$ and $p,\ord_A$ are even,}\\
\!0 &\text{in other case.} \endcases
$$
The result follows easily from this fact.\qed
\enddemo

\proclaim{Theorem 3.10} Let $A$ be a finite order Frobenius $k$-algebra.
Assume that the characteristic of $k$ does not divide $\ord_A$ and that $k$
has a primitive $\ord_A$-th root of unity $w$. For each $n\ge 1$ we write
$n = qe_A+s$, with $0\le s< e_A$. Then,

\smallskip

\item{1)} if $\ch(k)\ne 2$ and $\ord_A$ is odd,
$$
\align
\HH^n(TA) & = \HH_n(A)^* \oplus\bigoplus_{i=0}^q
\left(\HH^{n-ie_A}(A)\oplus \HH^{n-ie_A-1}(A)\right) \\
& \oplus\bigoplus_{j=2\atop j\text{ odd}}^{e_A}\bigoplus_{i=0}^{q-1}\left(
\H_{(j)}^{n-ie_A-j+1,0}(A)\oplus \H_{(j)}^{n-ie_A-j,0}(A)\right)\\
& \oplus\bigoplus_{j=2\atop j\text{ odd}}^{s+1}
\left(\H_{(j)}^{n-qe_A-j+1,0}(A)\oplus \H_{(j)}^{n-qe_A-j,0}(A) \right),
\endalign
$$
where $\HH^{-1}(A) = \H_{(s+1)}^{-1,0}(A) = 0$.

\smallskip

\item{2)} if $\ch(k)\ne 2$ and $\ord_A$ is even,
$$
\align
\HH^n(TA) & = \HH_n(A)^* \oplus\bigoplus_{i=0}^q
\left(\HH^{n-ie_A}(A)\oplus \HH^{n-ie_A-1}(A)\right) \\
& \oplus\bigoplus_{j=2\atop j\text{ odd}}^{e_A}\bigoplus_{i=0}^{q-1}\left(
\H_{(j)}^{n-ie_A-j+1,0}(A)\oplus \H_{(j)}^{n-ie_A-j,0}(A)\right)\\
&\oplus\bigoplus_{j=2\atop j\text{ odd}}^{s+1}\left( \H_{(j)}^{
n-qe_A-j+1,0} (A) \oplus \H_{(j)}^{n-qe_A-j,0}(A) \right)\\
& \oplus\bigoplus_{j=2\atop j\text{ even}}^{e_A}\bigoplus_{i=0}^{q-1}
\left( \H_{(j)}^{n-ie_A-j+1,\frac{e_A}{2}}(A)\oplus \H_{(j)}^{n-ie_A-j,
\frac{e_A}{2}} (A)\right)\\
& \oplus\bigoplus_{j=2\atop j\text{ even}}^{s+1} \left(
\H_{(j)}^{n-qe_A-j+1, \frac{e_A}{2}}(A)\oplus \H_{(j)}^{n-qe_A-j,
\frac{e_A}{2}}(A) \right),
\endalign
$$
where $\HH^{-1}(A) = \H_{(s+1)}^{-1,0}(A) =
\H_{(s+1)}^{-1,\frac{e_A}{2}}(A) = 0$.

\smallskip

\item{3)} if the characteristic of $k$ is $2$, then
$$
\align
\HH^n(TA) & = \HH_n(A)^* \oplus\bigoplus_{i=0}^q
\left(\HH^{n-ie_A}(A)\oplus \HH^{n-ie_A-1}(A)\right) \\
& \oplus\bigoplus_{j=2}^{e_A}\bigoplus_{i=0}^{q-1}\left(
\H_{(j)}^{n-ie_A-j+1,0}(A)\oplus \H_{(j)}^{n-ie_A-j,0}(A)\right)\\
& \oplus\bigoplus_{j=2}^{s+1} \left(\H_{(j)}^{n-qe_A-j+1,0}(A)\oplus
\H_{(j)}^{n-qe_A-j,0}(A) \right),
\endalign
$$
where $\HH^{-1}(A) = \H_{(s+1)}^{-1,0}(A) = 0$.
\endproclaim

\demo{Proof} This follows immediately from Theorem~3.8 and
Proposition~3.9.\qed
\enddemo

\remark{Remark 3.11} If $k$ does not have a primitive $\ord_A$-th root of
unity, we can apply the above theorem, using that $\ov{k}\ot \HH^*(T(A))
= \HH_{\ov{k}}^*(T(\ov{k}\ot A))$, where $\ov{k}$ is a suitable extension
of $k$ and $\HH_{\ov{k}}^*(T(\ov{k}\ot A))$ denotes the Hochschild
cohomology of $T(\ov{k}\ot A) = \ov{k}\ot T(A)$ as a $\ov{k}$-algebra.
\endremark

\remark{Remark 3.12} As it is well known, every finite dimensional Hopf
algebra $A$ is Frobenius, being a Frobenius homomorphism any right integral
$\varphi\in A^*\setminus\{0\}$. Moreover, by \cite{S, Proposition~3.6}, the
compositional inverse of the Nakayama map $\rho$ with respect to $\varphi$,
is given by
$$
\rho^{-1}(h) = \al(h_{(1)})\ov{S}^2(h_{(2)}),
$$
where $\al\in A^*$ is the modular element of $A^*$ and $\ov{S}$ is the
compositional inverse of $S$ (note that the automorphism of Nakayama
considered in \cite{S} is the compositional inverse of the considered by
us). Using this formula and that $\al\circ S^2 = \al$ it is easy to check
that $\rho(h) = \al(S(h_{(1)}))S^2(h_{(2)})$, and more generality, that
$$
\rho^l(h) = \al^{*l}(S(h_{(1)}))S^{2l}(h_{(2)}),
$$
where $\al^{*l}$ denotes the $l$-fold convolution product of $\al$. Since
$\al$ has finite order respect to the convolution product and, by the Radford
formula for $S^4$ (see \cite{S, Theorem~3.8}), the antipode $S$ has finite
order respect to the composition, $A$ has finite order. So, Theorem~3.10
applies to finite dimensional Hopf algebras.
\endremark

\proclaim{Theorem 3.13} Let $A$ be a finite dimensional Hopf algebra. Assume
that $A$ and $A^*$ are unimodular. It is hold that

\smallskip

\item{1)} if $\ch(k)\ne 2$ and $\ord_A = 1$,
$$
\HH^n(TA) = \HH_n(A)^* \oplus\bigoplus_{i=0}^n \HH^i(A),
$$

\smallskip

\item{2)} if the characteristic of $k$ is $2$ and $\ord_A = 1$, then
$$
\HH^n(TA)  = \HH_n(A)^* \oplus\bigoplus_{i=0}^n \HH^i(A)
\oplus\bigoplus_{i=0}^{n-1} \HH^i(A),
$$

\smallskip

\item{3)} if $\ch(k)\ne 2$ and $\ord_A = 2$,
$$
\HH^n(TA) = \HH_n(A)^* \oplus\bigoplus_{i=0}^n \HH^i(A) \oplus
\bigoplus_{i=0}^{n-1}\H_{(2)}^{i,1}(A).
$$
\endproclaim

\demo{Proof} By \cite{S, Corollary~3.20}, if $A$ and $A^*$ are unimodular,
then $\al$ is the counity $\ep$ of $A$ and $S^4 = id$. Hence, in this case,
$\rho = S^2$ and $\rho^2 = id$. The result follows immediately from this fact
and Theorem~3.10.\qed
\enddemo

\proclaim{Corollary 3.14} Let $A$ be a finite dimensional Hopf algebra. If
$A$ is semisimple and $\ch(k) = 0$, then
$$
\HH^n(TA) = \cases \HH_0(A)^* \oplus \HH^0(A) &\text{if $n = 0$,}\\
\HH^0(A) &\text{if $n > 0$.} \endcases
$$
\endproclaim

\demo{Proof} Since for each semisimple finite dimensional Hopf algebra over
a characteristic zero field $A$ it is hold that $A$ and $A^*$ are
unimodular Hopf algebra, the antipode $S$ of $A$ is involutive, the result
follows immediately from item~1) of Theorem~3.13 and the fact that and $A$
is separable.\qed
\enddemo

\proclaim{Theorem 3.15} Let $A$ be a  finite order Frobenius algebra.
Assume that the characteristic of $k$ does not divide $\ord_A$ and that $k$
has a primitive $\ord_A$-th root of unity $w$. Then $\HH^n(A) =
\HH^{n,0}_{(1)}(A)$ for all $n\ge 0$.
\endproclaim

\demo{Proof} By Propositions~3.4 and 3.9
$$
\align
H^n(Y_{(1)}^*) & = \cases \HH^n(A) & \text{if $n=0$,}\\ \HH^n(A) \oplus
\HH^{n-1}(A) & \text{if $n>0$.}\endcases \\
& = \cases \HH^{n,0}_{(1)}(A) & \text{if $n=0$,}\\ \HH^{n,0}_{(1)}(A)
\oplus \HH^{n-1,0}_{(1)}(A) & \text{if $n>0$.}\endcases
\endalign
$$
From this follows easily that $\HH^n(A) = \HH^{n,0}_{(1)}(A)$ for all $n\ge
0$, as desired.\qed
\enddemo

\example{Example 3.16} Let $k$ a field and $N$ a natural number. Assume
that $k$ has a primitive $N$-th root of unity $w$. Let $A$ be the Taft
algebra of order $N$. That is, $A$ is the algebra generated over $k$ by two
elements $g$ and $x$ subject to the relations $g^N = 1$, $x^N = 0$ and $xg
= wgx$. The Taft algebra $A$ is a Hopf algebra with comultiplication $\De$,
counity $\ep$ and antipode $S$ given by
$$
\alignat2
& \De(g) = g\ot g, &&\qquad  \De(x) = 1\ot x + x\ot g, \\
& \ep(g) = 1, && \qquad \ep(x) = 0, \\
& S(g) = g^{-1}, &&\qquad  S(x) = - xg^{-1}.
\endalignat
$$
Using that $t = \sum_{j=0}^{N-1} w^jg^jx^{N-1}$ is a right integral of $A$,
it is easy to see that the modular element $\al\in A^*$ verifies $\al(g) =
w^{-1}$ and $\al(x) = 0$. By Remark~3.14 the Nakayama map $\rho\:A\to A$ is
given by $\rho(g) = wg$ and $\rho(x) = w^{-1}x$. Hence, $A = A_0\oplus
\cdots\oplus A_{N-1}$, where
$$
\align
A_i & = \{a\in A:\rho(a) = w^{-i}a\}\\
& = \langle x^i,x^{i+1}g,\dots, x^{N-1}g^{N-i-1},g^{N-i},xg^{N-i+1},\dots,
x^{i-1}g^{N-1}\rangle.
\endalign
$$
Let $C_N = \{1,t,\dots,t^{N-1}\}$ be the cyclic of order $N$. It is easy to
see that $C_N$ acts on $A_0$ via $t\cdot x^i g^i = w^{-i}x^i g^i$ and $A$
is isomorphic to the skew product of $A_0\# C_N$. By Theorem~3.15 we have
$$
\HH^n(A) = \HH^{n,0}_{(1)}(A) = H^n(Y_{(1),0}^{0,*}),
$$
where $Y_{(1),0}^{0,*}$ is the complex introduced in the proof of
Theorem~3.9. Let us consider the action of $C_N$ on $\HH^q(A_0)$ is induced
by the action of $C_N$ on $\Hom_k(A_0^{\ot q},A_0)$ given by
$$
t \cdot \varphi(x^{i_1}g^{i_1}\ot \cdots\ot x^{i_q}g^{i_q}) = g^{N-1}
\varphi(t\cdot x^{i_1}g^{i_1}\ot \cdots\ot t\cdot x^{i_q}g^{i_q}) g.
$$
By \cite{G-G2, 3.2.7} there is a converging spectral sequence
$$
E^2_{pq} = \H^p(C_N,\HH^q(A_0)) \Rightarrow H^{p+q}(Y_{(1),0}^{0,*}).
$$
Now, $N$ is invertible in $k$, since $k$ has a primitive $N$-th root of
unity. Thus, the above spectral sequence collapses and
$$
\HH^n(A) = \H^0(C_N,\HH^n(A_0)) = \HH^n(A_0)^{C_N}.
$$
\endexample

\enddocument